\documentclass{article}

\usepackage{latexsym}
\usepackage{amsmath}
\usepackage{amssymb}
\usepackage{graphics}
\usepackage{amsthm}
\usepackage{epsfig}

\renewcommand{\Bbb}[1]{\mathbb{#1}}

\theoremstyle{plain}

\newtheorem{lem}{Lemma}
\newtheorem{cor}{Corollary}

\newtheorem{defn}{Definition}                  
\newtheorem{thm}{Theorem}

\def\MyDate{\the\month/\the\day/\the\year}

\title{The Casson-Walker-Lescop invariant and link invariants}
                                             \author{Jeff Johannes}                       
\date { }

\begin{document}

\begin{center}
 {\bf THE CASSON-WALKER-LESCOP INVARIANT AND LINK INVARIANTS}\\[1cm]

{\small

JEFF JOHANNES

{\em Department of Mathematical Sciences

University of Nevada Las Vegas

4505 Maryland Parkway

Box 454020

Las Vegas, NV  89154-4020

johannes@member.ams.org}\\[1in]

ABSTRACT}
\end{center}

{\small
\begin{quotation}
Formulas previously presented for the Casson-Walker invariant are generalized to Lescop's extension.  These formulas in terms of linking numbers and surgery coefficients compute the change in Lescop's invariant under crossing changes in a framed link presenting a 3-manifold.  This leads us to revisit an old formula for a coefficient of the Conway polynomial.  Finally we compute Lescop's invariant of several lens spaces and deduce some values for Dedekind sums.\\

\noindent {\em Keywords:}  3-manifolds, links, Casson-Walker-Lescop invariant, Conway polynomial, Dedekind sums.
\end{quotation}}

This version was compiled on \MyDate.

\section{Introduction}

In 1985 Casson introduced an invariant, $\lambda_c$, to distinguish homology 3-spheres by counting $SU(2)$ representations of their fundamental groups [AM].  In 1988 Walker carefully refined the details of Casson's work in order to extend Casson's definitions to an invariant, $\lambda_w$, of rational homology spheres [W].  Walker also established a method of defining $\lambda_w$ that relied only on the surgery presentation of the manifold and avoided the question of representations altogether.  In 1992 Lescop generalized the Casson-Walker invariant to the rational-valued invariant $\lambda$ of any 3-manifold [L].  Lescop's invariant is not interpreted in terms of representations and is only defined as a combinatorial invariant of a 3-manifold presented as surgery on a link.  

In [J] we presented a formula in terms of linking numbers and surgery coefficients for computing how the Casson-Walker invariant changes under crossing changes in framed links presenting 3-manifolds.  Theorem 1 of this current paper shows that after being adjusted appropriately, the Casson-Walker crossing change formula extends to Lescop's invariant of arbitrary 3-manifolds.  This allows for geometric computations of Lescop's invariant which bypass the use of the original page-long formulas in [L].  As an aside, we see a method of extracting a difference of the $z^{n+1}$ coefficients of the Conway polynomial of two $n$-component links that differ by a crossing change from information about Lescop's invariant of manifolds derived by particular surgeries on these links.  

As examples of the ease of computing Lescop's invariant using these new techniques, we compute the invariant for several lens spaces.  The two results presented in Section 3 are Theorem 2:  For any natural number $r$, \\
$\lambda(L(r^2+ 1, r)) = 0$, and Theorem 3:  For any natural numbers $n$, and $b$ 
\[\lambda(L(2n^2b^2 + 2nb + 1, 2nb^2)) = \frac{b^2 (n^3 - n)}{12}.\]

In both of these cases we computed Lescop's invariant without using the original combinatorial definition.  Hence we may set these values equal to the combinatorial presentation and solve for the most intricate piece, the Dedekind sums.  This technique produces the following two formulas for Dedekind sums:  Corollary 1:  For any natural number $r$,
$s(r^2+1, r) = \frac{ r^2 - 3r + 2}{12r}$  and Corollary 2:  For any natural numbers $n$, and $b$
\[s(2n^2b^2 + 2nb + 1, 2nb^2) = {\frac {2\,{n}^{2}{b}^{2}-3\,n{b}^{2}+1}{12n{b}^{2}}}. \]

The techniques used to compute Lescop's invariant in these examples can be used to compute the invariant of any 3-manifold presented as a two component link ${\bf L}$ with $(s, -s)$ framings on the components for any rational number $s$.  

\section{A Crossing Change Formula for the Casson-Walker-Lescop Invariant}

\subsection{Background}

Let us set some notation and recall some formulas.  For any abelian group $\Gamma$, $|\Gamma|$ is defined to be the order of the group and defined
to be
zero if $\Gamma$ is infinite.  Let ${\bf L} = (L_1, L_2, . . ., L_n)$ be an ordered $n$-component link in $S^3$, let ${\bf s} = (s_1, s_2, . . ., s_n)$ be a rational vector, let ${\bf L_s}$ be the 3-manifold obtained by ${\bf s}$
surgery on
${\bf L}$, and let $A({\bf L}, {\bf s})$ be the associated linking matrix.  Let $\lambda({\bf L}_{\bf s})$ be the 
Casson-Walker-Lescop invariant of ${\bf L_s}$, and if ${\bf L_s}$ is a
rational homology sphere, or equivalently if det$(A({\bf L}, {\bf s})) \ne 0$, let
$\lambda_w({\bf L}_{\bf s})$ be the Casson-Walker invariant [AM, W] of ${\bf L_s}$.

As is now standard, given an invariant of links, $\rho$, we may extend it to an invariant of links
with one transverse double point by defining \hbox{$\rho(\times) = \rho(+) - \rho(-)$,} that is by
taking
the difference between the positive and negative crossing resolutions of the double point [BL].  We
shall refer to this difference as a crossing change.  
For a given singular link with one transverse double point ${\bf L}$, we will refer to ${\bf
L}^+$ as the link
resulting from a positive resolution of the double point, and ${\bf L}^-$ as the negative
resolution.  So that the singular link has the same number of components as each of its resolutions, we limit
to internal crossing changes, that is, crossing changes within a component.  We define a singular link to be {\em disjoint} if none of the singularities arise as  intersections {\em between} components.  This limitation is
equivalent to allowing only link homotopy transitions from one link to another.  

In [J] we presented a crossing change formula for computing the difference of the Casson-Walker invariants of two 3-manifolds presented as framed links that differ by one crossing change within the first component of the link.  When we perform a
crossing change we can consider passing through an intermediate stage of a link with one singular point from a self-intersection in the first component.  If we have precisely one singularity there are two ``lobes'' determined.  To
see one lobe start at the singularity and traverse the knot until returning to the singularity.  The other lobe 
is the remaining portion of the component.  We will denote the two lobes of $L_1$ as $L_1^a$
and $L_1^b$.  The process of separating the lobes is called ``smoothing'', and it is done by
replacing $\left(\includegraphics{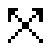}\right)$ with
$\left(\includegraphics{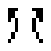}\right)$ locally.  Given
this
background the crossing change
formula can be restated.

For ${\bf L^-}$ and ${\bf L^+}$, two links in $S^3$ with
the same framing which differ by a crossing change in the first component, where 
det$(A({\bf L^-}, {\bf s})) = $ det$(A({\bf L^+}, {\bf s})) \ne 0$,
\[(*) \  \lambda_w({\bf L}_{\bf s}^-) - \lambda_w({\bf L}_{\bf s}^+) = \frac {2 \begin{vmatrix}
l & k_{12}^a & \dots & k_{1n}^a\\
k_{12}^b & s_2 & \dots & n_{2n}\\
\vdots & \vdots & \ddots & \vdots \\
k_{1n}^b & n_{2n} & \dots & s_n\\
\end{vmatrix} }
{\begin{vmatrix}
s_1 & n_{12} & \dots & n_{1n}\\
n_{12} & s_2 & \dots & n_{2n}\\
\vdots & \vdots & \ddots & \vdots \\
n_{1n} & n_{2n} & \dots & s_n\\
\end{vmatrix}}  \]
where $s_i$'s are the surgery coefficients, $n_{ij}$'s are the linking numbers between the
components, $k_{1j}^a$'s are the linking numbers between one lobe and the other components, and
$l$ is the linking number between the two smoothed lobes.

In 1992 Lescop [L] fully generalized the Casson-Walker invariant to the 
rational-valued invariant
$\lambda$ of
any 3-manifold.  If M is a rational homology sphere then
$\lambda(M) = \frac{|H_1(M)|}{2} \lambda_w(M)$ (which, if we restrict further to integral
homology spheres, provides that $\lambda = \lambda_c$, Casson's original invariant of integral homology spheres).  Lescop's invariant is no longer
interpreted in terms of representations and is now only defined as a combinatorial invariant of a 
3-manifold presented as surgery on a link.  



\subsection{Lescop's surgery formula}

In order to consider the crossing change formula in the context of Lescop's extension of the Casson-Walker invariant, we should recall some facts and formulas from [L].  
In the general case the surgery formula for Lescop's invariant is rather complicated.  For
our purposes, we will only use
the following formula for surgery in $S^3$.  (This entire subsection is taken
from [L, \S 1.7] and uses the notation used there.)

We start with a framed $n$ 
component link $({\bf L}, {\bf s})$ in $S^3$.  Let $N = \{1, \ldots,
n\}$.   Let $I$ be used for an index set that is a subset of $N$; let $\sharp I$ denote the number
of elements in $I$.  Given a matrix $M$ we define sign($M$) to be $(-1)^{b_-(M)}$ where $b_-
(M)$ is the number of negative eigenvalues of $M$.  The signature of a Hermitian matrix $M$ is
defined to
be $b_+(M) - b_-(M)$, where $b_+(M)$ is the number of positive eigenvalues.   
We will indicate
restricting a link to components in the index set $I$ by ${\bf L}_I$.  Similarly restricting a
matrix, that
is taking only the columns and rows with indices in $I$, is
indicated by $M_I$. 




There are two more intricate definitions we need.  

\begin{defn}
Let $q$ be an integer; let $p$ be an integer or a mod $q$ congruence class of an integer.  The {\em Dedekind sum} $s(p,q)$ is the following rational number:
\[s(p,q) := \sum_{i=1}^{|q|} \left(\left(\frac{i}{q}\right)\right) \left(\left(\frac{pi}{q}\right)\right)\]
where $((x)) =  \begin{cases}
0 & \text{if }    x \in \Bbb{Z} \\
x - [x] - \frac{1}{2} & otherwise
\end{cases}$
\end{defn}

We also must define a 
revision of the linking matrix, denoted $A(({\bf L}, {\bf s})_{N \backslash I}; I))$.  This matrix $A(({\bf L}, {\bf s})_{N \backslash I}; I))$ is a $(n - \sharp I) \times (n - \sharp I)$ matrix.    Let $n'_{ij}$
be the entries of this matrix where $i, j \in N \backslash I$ and
\[n'_{ij} =
\begin{cases}
n_{ij} &  i \ne j \\
s_{i} + \sum_{k \in I}n_{ki} & i = j
\end{cases}\] 
Here is the Conway polynomial formula for  Lescop's invariant:
\begin{multline*}
\lambda({\bf L_s}) = \text{sign}(A({\bf L}, {\bf s}))
\left(\prod_{i=1}^n q_i \right) 
\sum_{\{I | I \ne \emptyset, I \subseteq N\}}
\text{det}(A(({\bf L}, {\bf s})_{N \backslash I}; I)) a_1({\bf L}_I) \\
+ \text{sign}(A({\bf L}, {\bf s}))
\left(\prod_{i=1}^n q_i \right)
\sum_{\{I | I \ne \emptyset, I \subseteq N\}}
\frac{\text{det}(A(({\bf L}, {\bf s})_{N \backslash I}))(-1)^{\sharp I}\Theta(A(({\bf L}, {\bf
s})_I))}{24}\\
 + |H_1({\bf L_s})|\left(\frac{\text{signature}(A({\bf L}, {\bf s}))}{8} + \sum_{i=1}^n \frac{s(p_i,q_i)}{2}\right)
\end{multline*}
In this formula, $a_1$ is the coefficient on $z^{n+1}$ in the Conway polynomial, and
\[\Theta(A_I) = \begin{cases}
\Theta_b(A_I) + \frac{q_i^2 + 1}{q_i^2} & I = \{i\}\\
\Theta_b(A_I) - 2n_{ij} & I = \{i, j\}\\
\Theta_b(A_I) & \sharp I > 2
\end{cases}\]
where
\[\Theta_b(A_I) = \sum_\Xi
\text{Lk}_c(A_J)n_{ig(1)} n_{g(1)g(2)}\dots n_{g(\sharp (I \backslash J) - 1) g(\sharp (I
\backslash J))} n_{g(\sharp (I \backslash J)j}\]
in which $\Xi = {\{(J, i,j,g) | J \subseteq I, J \ne \emptyset, (i,j) \in J^2, g \in \sigma_{I
\backslash J}\}}$ and 
\[\text{Lk}_c(A_J) = \sum_{\sigma \in \sigma_J}
(n_{\sigma(1)\sigma(2)}n_{\sigma(2)\sigma(3)} \dots n_{\sigma(j-1)\sigma(j)}
n_{\sigma(j)\sigma(1)})\]

\subsection{Crossing Changes and the Casson-Walker-Lescop Invariant}
 
Multiplying equation $(*)$ in Subsection 2.1 by $\frac{|H_1(M)|}{2}$ produces a crossing change formula for Lescop's invariant restricted to rational homology spheres where both the Casson-Walker invariant and Lescop's extension are defined.  We next extend this
formula fully to Lescop's invariant of arbitrary 3-manifolds.        

\begin{thm} For all framed links in $S^3$, 
\[\lambda({\bf L}^-_{\bf s}) - \lambda({\bf L}^+_{\bf s}) = \text{sign}(A({\bf L}, {\bf s}))
\prod q_i \begin{vmatrix}
l & k_{12}^a & \dots & k_{1n}^a\\
 k_{12}^b & s_2 & \dots & n_{2n}\\
\vdots & \vdots & \ddots & \vdots \\
k_{1n}^b & n_{2n} & \dots & s_n\\
\end{vmatrix}  \]
\end{thm}

\begin{proof}
Approach: In the rational homology sphere case, we have the following crossing change formula $(*)$
for the Casson-Walker invariant from Subsection 2.1:
\[(*) \ \lambda_w({\bf L}^-_{\bf s}) - \lambda_w({\bf L}^+_{\bf s}) = \frac {2 \begin{vmatrix}
l & k_{12}^a & \dots & k_{1n}^a\\
k_{12}^b & s_2 & \dots & n_{2n}\\
\vdots & \vdots & \ddots & \vdots \\
k_{1n}^b & n_{2n} & \dots & s_n\\
\end{vmatrix} }
{\begin{vmatrix}
s_1 & n_{12} & \dots & n_{1n}\\
n_{12} & s_2 & \dots & n_{2n}\\
\vdots & \vdots & \ddots & \vdots \\
n_{1n} & n_{2n} & \dots & s_n\\
\end{vmatrix}}  \]
Recall that for rational homology spheres, $\lambda(M) = \frac{|H_1(M)|}{2} \lambda_w(M)$. 
Using the fact that $|H_1(M)| = \prod_i (q_i) $sign$(A({\bf L}, {\bf s}))$det$(A({\bf L}, {\bf s}))$, we may multiply the above equation $(*)$ by $\frac{|H_1(M)|}{2}$ to yield the crossing change formula for Lescop's
invariant in the case of
rational homology spheres, i.e. when det$(A({\bf L}, {\bf s})) \ne 0$:
\[\lambda({\bf L}^-_{\bf s}) - \lambda({\bf L}^+_{\bf s}) = \text{sign}(A({\bf L}, {\bf s}))
\prod q_i \begin{vmatrix}
l & k_{12}^a & \dots & k_{1n}^a\\
 k_{12}^b & s_2 & \dots & n_{2n}\\
\vdots & \vdots & \ddots & \vdots \\
k_{1n}^b & n_{2n} & \dots & s_n\\
\end{vmatrix}  \]  
Notice that if we divide $\lambda({\bf L}^-_{\bf s}) - \lambda({\bf L}^+_{\bf s})$ by sign$(A({\bf L}, {\bf s}))\prod q_i$ that the result is a polynomial in the components of ${\bf s}$, for ${\bf s}$ such that det$(A({\bf L}, {\bf s})) \ne 0$.

If we can show that the general crossing change formula for Lescop's invariant divided by
sign$(A({\bf L}, {\bf s}))\prod q_i$ must be
a polynomial in ${\bf s}$, then the two polynomials (the one obtained by rescaling the Casson-Walker
invariant as above and the one
directly from Lescop's formula)
will agree on the set det$(A({\bf L}, {\bf s})) \ne 0$, which is dense in $\Bbb{Q}^n$. 
Therefore these two
polynomials will agree everywhere in $\Bbb{Q}^n$ (i.e. their difference is zero).  So, the crossing
change formula in the rational homology sphere case will extend to arbitrary 3-manifolds. 

Taking this approach, it now suffices to show that the general crossing change formula for Lescop's
invariant divided by
sign$(A({\bf L}, {\bf s}))\prod q_i$ is
polynomial for every framed link.  We will do this via Lescop's formula in terms of the Conway
polynomial as presented in Subsection 2.2. The goal is to find a formula for $\lambda({\bf L}^-_{\bf s}) - \lambda({\bf
L}^+_{\bf
s})$.  
Recall
\begin{multline*}
\lambda({\bf L_s}) = \text{sign}(A({\bf L}, {\bf s}))
\left(\prod_{i=1}^n q_i \right) 
\sum_{\{I | I \ne \emptyset, I \subseteq N\}}
\text{det}(A(({\bf L}, {\bf s})_{N \backslash I}; I)) a_1({\bf L}_I) \\
+ \text{sign}(A({\bf L}, {\bf s}))
\left(\prod_{i=1}^n q_i \right) 
\sum_{\{I | I \ne \emptyset, I \subseteq N\}}
\frac{\text{det}(A(({\bf L}, {\bf s})_{N \backslash I}))(-1)^{\sharp I}\Theta(A(({\bf L}, {\bf
s})_I))}{24}\\
 + |H_1({\bf L_s})|\left(\frac{\text{signature}(A({\bf L}, {\bf s}))}{8} + \sum_{i=1}^n \frac{s(p_i,q_i)}{2}\right)
\end{multline*}
If we compute $\lambda({\bf L}_{\bf s}^-) - \lambda({\bf L}_{\bf s}^+)$ from this formula for Lescop's invariant, because ${\bf L}^+$ and ${\bf L}^-$ are in the same link homology class (i.e. have all the same linking numbers) the terms involving only linking numbers cancel.  This leaves us with
\begin{multline*}
\lambda({\bf L}_{\bf s}^-) - \lambda({\bf L}_{\bf s}^+) = \\
 \text{sign}(A({\bf L}, {\bf s}))
\left(\prod_{i=1}^n q_i \right)
\sum_{\{I | I \ne \emptyset, I \subseteq N\}}
\text{det}(A(({\bf L}, {\bf s})_{N \backslash I}; I)) [a_1({\bf L}_I^-) - a_1({\bf L}_I^+)] 
\end{multline*}
If we consider instead $\frac{\lambda({\bf L}_{\bf s}^-) - \lambda({\bf L}_{\bf s}^+)}{\text{sign}(A({\bf L}, {\bf s}))
\left(\prod_{i=1}^n q_i \right)}$, we see, as desired, that this quantity is a polynomial in the variables $s_i$ just as the
other
crossing change formula. 
\end{proof}

Note:  what is proven in this argument is useful to note in a general context.  Any polynomial formula that holds for $\left(\lambda_w({\bf L}_{\bf s}^1) - \lambda_w({\bf L}_{\bf s}^2)\right)(\text{det}(A({\bf L}, {\bf s}))$ for all rational homology spheres ${\bf L}_{\bf s}^1$ and ${\bf L}_{\bf s}^2$ in the same link homology class can be extended to a formula that holds for $\frac{\lambda({\bf L}_{\bf s}^1) - \lambda({\bf L}_{\bf s}^2)}{\text{sign}(A({\bf L}, {\bf s}))
\left(\prod_{i=1}^n q_i \right)}$ for arbitrary 3-manifolds.  

\subsection{Connections to Coefficients of the Conway Polynomial}

It is perhaps interesting to note that via Theorem 1 we may extract from the Casson-Walker-Lescop invariant a difference of the $a_1$ coefficients of the Conway polynomial for the presenting links.  If we consider setting $s_i = -\sum_{k \in N, k \ne i} n_{ik}$, then for any $I \ne N, A(({\bf L}, {\bf s})_{N \backslash I}; I)$ has the property that the sum of each row is zero.  Any such matrix has determinant zero (consider that eigenvector $(1, 1,  \dots 1)$ has eigenvalue 0).  Therefore for this value of ${\bf s}$ each of the terms in 
\[\frac{\lambda({\bf L}_{\bf s}^-) - \lambda({\bf L}_{\bf s}^+)}{\text{sign}(A({\bf L}, {\bf s}))
\left(\prod_{i=1}^n q_i \right)} = \sum_{\{I | I \ne \emptyset, I \subseteq N\}}
\text{det}(A(({\bf L}, {\bf s})_{N \backslash I}; I)) [a_1({\bf L}_I^-) - a_1({\bf L}_I^+)] \]
is zero except for $I = N$.  Therefore we have 
\[\frac{\lambda({\bf L}_{\bf s}^-) - \lambda({\bf L}_{\bf s}^+)} {\text{sign}(A({\bf L}, {\bf s}))
\left(\prod_{i=1}^n q_i \right)} = a_1({\bf L}^-) - a_1({\bf L}^+).\]

Combining this with Theorem 1 yields 
\[ a_1({\bf L}^-) - a_1({\bf L}^+) = \begin{vmatrix}
l & k_{12}^a & \dots & k_{1n}^a\\
 k_{12}^b & s_2 & \dots & n_{2n}\\
\vdots & \vdots & \ddots & \vdots \\
k_{1n}^b & n_{2n} & \dots & s_n\\
\end{vmatrix}\]
This is an obvious consequence of the formula for $a_0$ presented in [H] produced from an entirely different analysis.

\section{Examples and Consequences}

Let us examine some examples of applications of the formulas in Section 2.  We will focus on lens spaces because they are a class of spaces with several accessible presentations.  For our purposes we will restrict to writing our lens spaces $L(p,q)$ with positive $p$ and $q$.  The simplest surgery presentation of $L(p, q)$ is as $\frac{p}{q}$ surgery on the unknot.  
\[\includegraphics{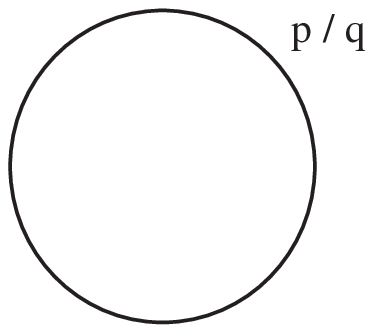}\]
From this presentation we may use Lescop's formula 1.4.8 [L] to compute \[\lambda(L(p,q)) = q \left( -\frac{1}{24} - \frac{p^2 + 1}{24q^2} \right) + \frac{p}{8} + \frac{p  \ s(p, q)}{2}\]  
Computing this in general requires computing Dedekind sums.  In special cases we may avoid this and as a consequence infer some formulas for Dedekind sums.  

Let us first examine a very special case.  We may present $L(r^2 + 1, r)$ as surgery on the Hopf link with $r$ and $-r$ coefficients.  
\[\includegraphics{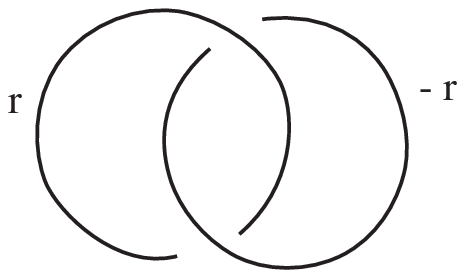} \]
Because this framed link is isotopic to the mirror image with the components interchanged, which here is the same as being isotopic to the mirror image with negative surgery coefficients, we find that $\lambda(L(r^2+ 1, r)) = - \lambda(L(r^2+ 1, r))$, and therefore 

\begin{thm}  For any natural number $r$,
$\lambda(L(r^2+ 1, r)) = 0$.  
\end{thm}

From this we can solve the above formula for $\lambda(L(p,q))$ for $s(p,q)$ to get the following:

\begin{cor}  For any natural number $r$,
\[s(r^2+1, r) = \frac{ r^2 - 3r + 2}{12r}.\]
\end{cor}

Applying this same philosophy and using Theorem 1 we can produce more interesting computations.  (Some of the ideas in the following example were previously worked out in [KL] for the Casson-Walker invariant with $s = 0$.)  Let us  generalize the Hopf link to the following link, denoted $T(n)$, with rational surgery coefficients $s = \frac{a}{b}$ and $-s$:
\[\includegraphics{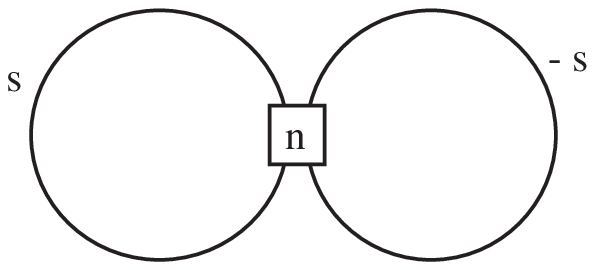}\]
(Here the $n$ in the box indicates that the two components have linking number $n$.)  This link is obviously link  homotopic to its mirror image $T'(n)$(as all two component links are) and interchangeable.  Therefore, we may find a link homotopy path to the manifold with reversed orientation.   
$T(n)$ is isotopic to the following presentation (here for $n = 4$):
\[\includegraphics{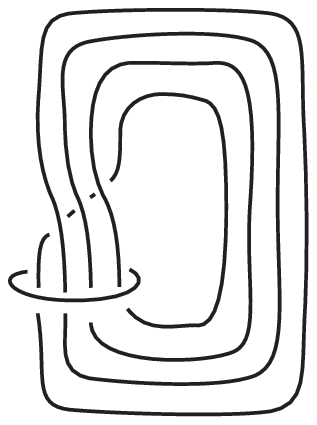}\]
From this presentation we may see the link homotopy to the mirror image is obtained by changing all of the self-crossings and twisting the circle component up.  Hence we may compute $\lambda(T(n), (s, -s)) - \lambda(T'(n), (-s, s))$.  Applying Theorem 1 to the $n-1$ crossing changes to convert $T(n)$ to $T'(n)$ we find 

\begin{multline*}
\lambda(T(n), (s, -s)) - \lambda(T'(n), (-s, s)) =  \sum_{i = 1}^{n - 1} \text{sign} \left( \begin{vmatrix}
s & n\\
n & -s\\
\end{vmatrix} \right)
b^2 \begin{vmatrix}
0 & n - i\\
i & s\\
\end{vmatrix}\\
 = - b^2 \sum_{i = 1}^{n - 1} (i^2 - ni) = \frac{b^2 (n^3 - n)}{6}
\end{multline*}
Therefore, $\lambda(T(n), (s, -s)) = \frac{b^2 (n^3 - n)}{12}$.  So that we can view this link surgery as gluing two tori togeher to form a lens space, we require that the surgery map is given by a homeomorphism.  The surgery map for the first component induces a map on first homology given by the matrix 
$\begin{bmatrix}
a & n \\
b & 1 
\end{bmatrix}$ (after choosing orientations appropriately).  The determinant of this matrix is then $a - nb$.  Setting this equal to one to give an isomorphism yields the condition $a = nb + 1$.  Therefore if $a = nb + 1$ then $(T(n), (s, -s))$ presents a lens space.  Which lens space this produces can be seen by a sequence of Kirby calculus moves [K], [R2].  Before performing this Kirby calculus manipulation, we will need two lemmas.

\begin{lem}
The following surgery presentations are equivalent for any rational framing $s$.
\[\epsfig{file=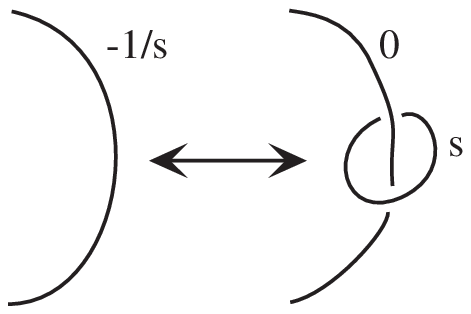,height=2cm}\]
\end{lem}

\begin{proof}
This result is similar to but slightly different from several standard surgery modifications presented in the literature.  A proof is presented here for completeness.  

Verifying this lemma is merely a matter of carefully tracing through surgery definitions. Let $s = \frac{a}{b}$.  Consider a ``bulging neighborhood" of the zero-framed component to include the $s$ framed component.  If we perform the two indicated surgeries within this solid torus, we are returned with a solid torus.  The remaining question is how this torus is reattached to the exterior of this neighborhood.  In other words, we see that the surgery instructions can be reduced to one local component as indicated on the left, but we have yet to determine the framing of that component.  

\[\epsfig{file=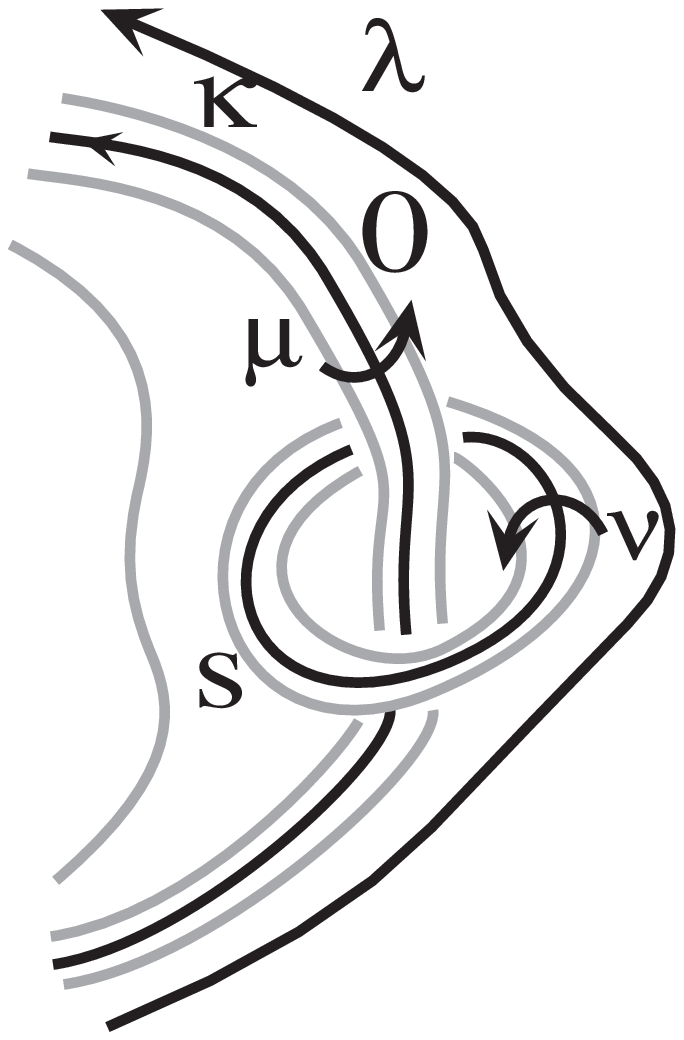,height=4cm}\]

For any surgeries on this solid torus, the resulting first homology will be generated by the meridians of the two link components, denote them $\mu$ and $\nu$, as indicated, and the longitude of the entire torus, denoted $\lambda$.  To help keep track of curves, let us also denote by $\kappa$ the longitude of the 0-framed component.  Because the $s$-framed component links the other component trivially once, the longitude of the $s$-framed component is isotopic to $\mu$, the meridian of the 0-framed component. 

Performing the 0-surgery reveals the relation $\kappa = 0$ in first homology, while performing the $s$-surgery provides $b\mu + a\nu = 0$.   It is apparent from the diagram that $\lambda = \kappa - \nu$.  This reduces to $\lambda = -\nu$.  Substituting into the first relation gives us $b\mu - a\lambda = 0$.  In fact, $\mu$ and $\lambda$ are isotopic to the meridian and longitude of the solid torus, so we recognize this as $\frac{b}{-a} = -\frac{1}{s}$ surgery on the core of that torus, as indicated in the surgery diagram on the left.
\end{proof}

\begin{lem}
The following surgery presentations are equivalent for any rational framings $s$, $t$, and $u$.  
\[\epsfig{file=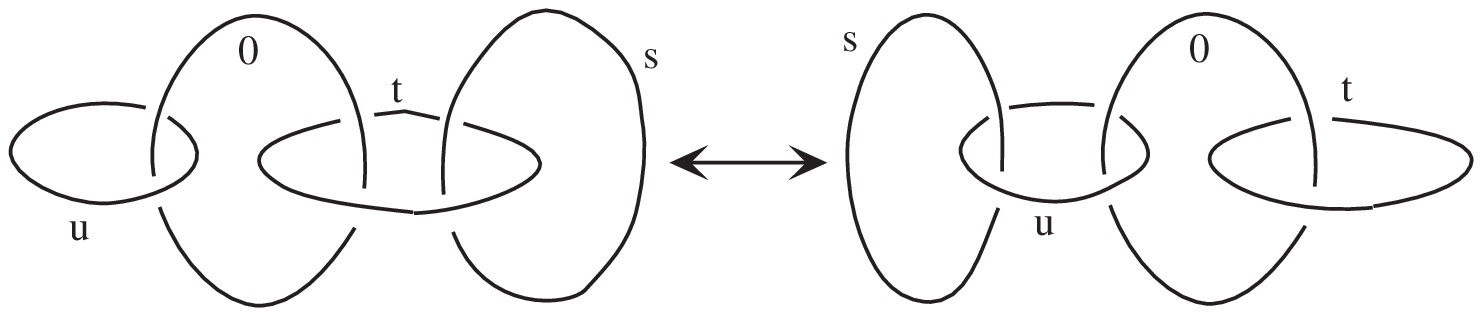,width=4.5in}\]
\end{lem}
\begin{proof}
Start with the presentation on the left.  Apply a Kirby move of the second kind, sliding the component with $s$ framing over the component with $0$ framing.  This unlinks the $s$ component from the $t$ component, but links the $s$ component with the $u$ component.  Because the framing of the central component is $0$ the framing of the $s$ component does not change in the process.
\end{proof}

Now for the Kirby calculus manipulations.  Let us set this result as another lemma.

\begin{lem}
The manifold presented by $(T(n), s, -s)$ is lens space \\ $L(2n^2b^2 + 2nb + 1, 2nb^2)$.
\end{lem}
\begin{proof}
First write out the surgery coefficients explicitly so that we may see them clearly.  
\[\includegraphics{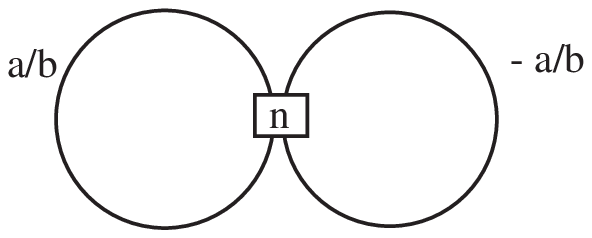}\]
Then insert a trivial component.
\[\includegraphics{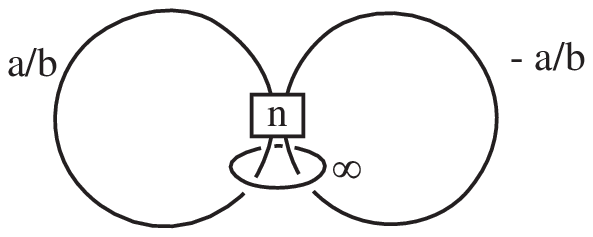}\]
Perform $n$ left-handed twists to separate the link.
\[\includegraphics{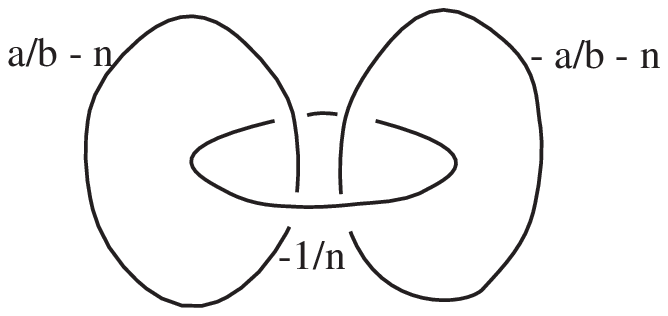}\]
Next apply Lemma 1 to the central component, changing the framing of the central component to zero and introducing a new component. 
\[\includegraphics{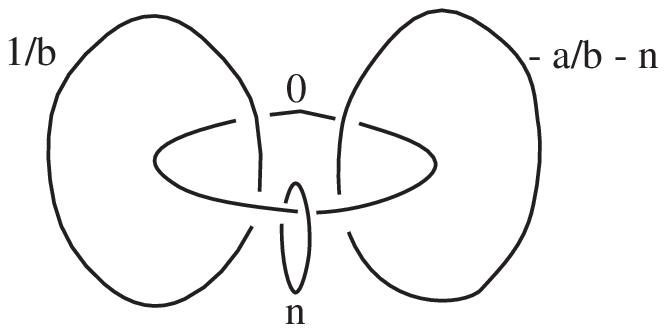}\]
Perform the same maneuver to the left component.
\[\includegraphics{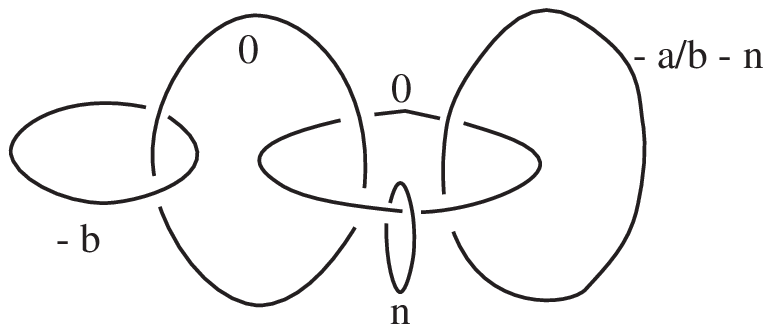}\]
Slide the rightmost component over the leftmost zero-framed component (by Lemma 2) and link it to the now $-b$ component.  
\[\includegraphics{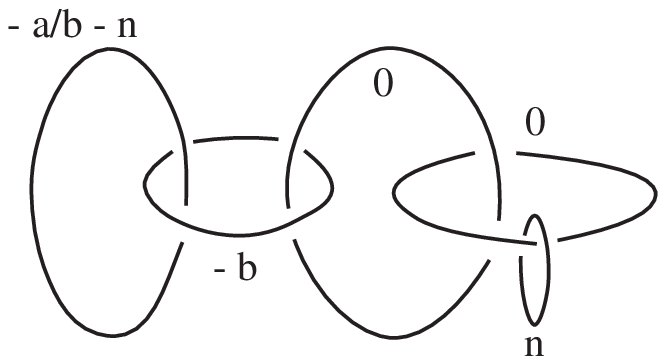}\]
At this point we now collapse the diagram by applying Lemma 1 in the opposite direction.  Eliminate one zero framed component 
\[\includegraphics{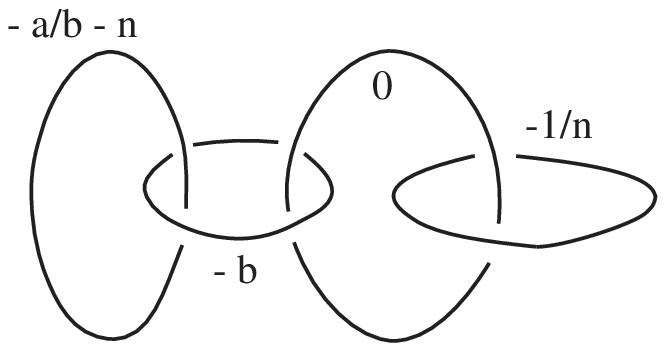}\]
then another.
\[\includegraphics{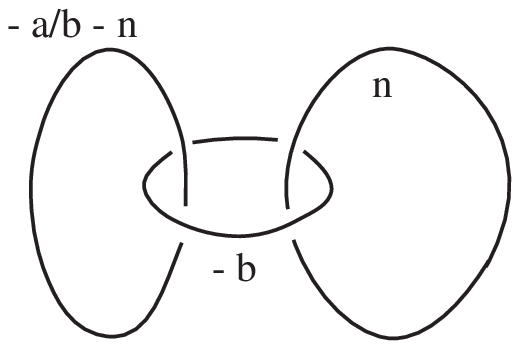}\]

It is a well-known result (stated in [R1] and [PS]) that if a manifold is presented as a ``chain" with integer surgery coefficients (or a rational coefficient on one end) then the manifold is a lens space.  In particular, for $a_i$ and $b_{n+1}$ integers, the following manifold 
\[\includegraphics{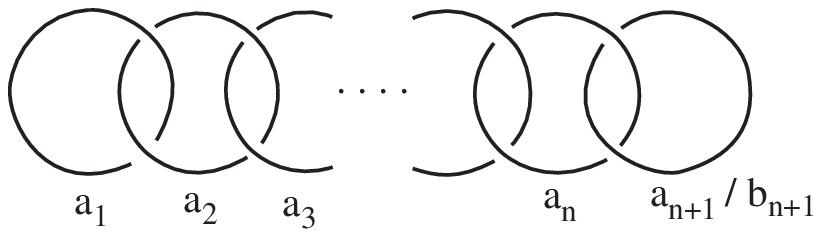}\]
is a presentation of $L(p,q)$ where $p$ and $q$ are given by the continued fraction
\[\frac{p}{q} = a_1 - \frac{1}{a_2 - \frac{1}{a_3 - \dots - \frac{1}{a_n - \frac{b_{n+1}}{a_{n+1}}}}}\]

By this result, the lens space we have here is given by 
\[\frac{p}{q} = n - \frac{1}{-b + \frac{b}{a + nb}}\]
Recalling the condition that $a = nb + 1$ this gives us 
\[\frac{p}{q} = \frac{2n^2b^2 + 2nb + 1}{2nb^2}\]

Therefore, $(T(n), s, -s) = L(2n^2b^2 + 2nb + 1, 2nb^2)$.
\end{proof}

We therefore have the following:

\begin{thm} For any natural numbers $n$, and $b$ 
\[\lambda(L(2n^2b^2 + 2nb + 1, 2nb^2)) = \frac{b^2 (n^3 - n)}{12}.\]
\end{thm}
Using this with Lescop's formula we find another formula for some Dedekind sums:

\begin{cor} For any natural numbers $n$, and $b$
\[s(2n^2b^2 + 2nb + 1, 2nb^2) = {\frac {2\,{n}^{2}{b}^{2}-3\,n{b}^{2}+1}{12n{b}^{2}}}. \]
\end{cor}

Applying these same techniques we can compute $\lambda({\bf L}_{(s, -s)})$ for any two component link ${\bf L}$ and any rational number $s$.  Because two component link homotopy classes are determined by linking numbers there is a homotopy between any link and its interchanged mirror image.  Using this homotopy we can compute $\lambda({\bf L}_{(s, -s)}) - \lambda(\bar{{\bf L}}_{(-s, s)}) = 2 \lambda({\bf L}_{(s, -s)})$ as we did in the above examples.  

\section*{Acknowledgements}

Some of the original work on this paper was completed while I was a graduate student at Indiana University.  I would like to thank my advisor Charles Livingston who encouraged me to follow many of the directions pursued in this paper.

\end{document}